%% file: main.tex
\documentclass{ifacconf}

\usepackage{graphicx}      
\usepackage{natbib}        
\usepackage{amsmath}
\usepackage{amssymb}
\usepackage{latexsym}
\usepackage{mathrsfs}
\usepackage{amscd}
\usepackage{sidecap}
\usepackage{graphicx}
\usepackage{color}
\usepackage{stmaryrd}
\usepackage{algorithm,algpseudocode}
\graphicspath{{pics/}}
\def\U{ {\mathcal {U}}}

\newcommand{\R}{\mathbb R}
\newcommand{\J}{\mathcal J}

\def\be#1\ee{\begin{equation}#1\end{equation}}

\newcommand{\bx}{\mathbf{x}}
\newcommand{\bv}{\mathbf{v}}
\newcommand{\bu}{\mathbf{u}}

\renewcommand{\theequation}{\arabic{section}.\arabic{equation}}
\renewcommand{\thetable}{\arabic{section}.\arabic{table}}
\renewcommand{\thefigure}{\arabic{section}.\arabic{figure}}

\def\S{\sigma}

\newcommand{\bd}{\begin{displaymath}}
\newcommand{\ed}{\end{displaymath}}
\newcommand{\ba}{\begin{eqnarray}}
\newcommand{\ea}{\end{eqnarray}}

\def\R{\mathbb{R}}

\def\R{\mathbb{R}}

\newcommand{\prt}[1]{\left(#1\right)}
\newcommand{\brk}[1]{\left[#1\right]}
\newcommand{\norm}[1]{\left\|#1\right\|}
\newcommand{\normsq}[1]{\left\|#1\right\|^2}
\newcommand{\inner}[2]{\left\langle#1,#2\right\rangle}

\usepackage{graphicx}      
\usepackage{natbib}        
\usepackage{amsmath}
\begin{document}
\begin{frontmatter}

\title{Optimal consensus control of the Cucker-Smale model} 

\author[First]{Rafael Bailo}
\author[Second]{Mattia Bongini}
\author[First]{Jos\'e A. Carrillo} 
\author[First]{Dante Kalise}

\address[First]{Department of Mathematics, Imperial College London, South Kensington Campus,  London SW7  2AZ, United Kindgom (e-mail: \{r.bailo,carrillo,dkaliseb\}@ic.ac.uk)}
\address[Second]{CEREMADE, Universit\'e Paris Dauphine, Place du Mar\'echal de Lattre de Tassigny, 75775 Paris Cedex 16, France (email: bongini@ceremade.dauphine.fr)}


\begin{abstract}                
We study the numerical realisation of optimal consensus control laws for agent-based models. For a nonlinear multi-agent system of Cucker-Smale type, consensus control is cast as a dynamic optimisation problem  for which we derive first-order necessary optimality conditions. In the case of a smooth penalisation of the control energy, the optimality system is numerically approximated via a gradient-descent method. For sparsity promoting, non-smooth $\ell_1$-norm control penalisations, the optimal controllers are realised by means of heuristic methods. For an increasing number of agents, we discuss the approximation of the consensus control problem by following a mean-field modelling approach.
\end{abstract}


\end{frontmatter}
\section{Introduction}

This paper addresses centralised control problems for second-order, nonlinear, multi-agent systems (MAS). In such dynamics, the state of each agent is characterised by a pair $(x_i,v_i)$, representing variables which we refer to as \textsl{position} and \textsl{velocity}, respectively. The uncontrolled system consists of simple binary interaction rules between the agents, such as attraction, repulsion, and alignment forces. This commonly leads to self-organisation phenomena, flocking or formation arrays. However, this behaviour strongly depends on the cohesion of the initial configuration of the system, and therefore control design is relevant in order to generate an external intervention able to steer the dynamics towards a desired configuration.  For second-order MAS, it is of particular interest the study of consensus emergence and control. In this context, we understand consensus as a travelling formation in which every agent has the same velocity. Self-organised consensus emergence for the Cucker-Smale model, see \cite{CS}, has been already characterised in \cite{HaLiu,CFRT}. The problem of consensus control for Cucker-Smale type models has been discussed in \cite{CFPT,bfk}. A related problem is the design of controllers achieving a given formation, which has been previously addressed in \cite{perea,borzi}. 

\textsl{Contributions.} In this work we focus on the design of centralised control laws enforcing consensus emergence. For this, we cast the consensus control problem in the framework of optimal control theory, for which an ad-hoc computational methodology is presented. We consider a finite horizon control problem, in which the deviation of the population with respect to consensus is penalised along a quadratic control term. We derive first-order optimality conditions, which are then numerically realised via the Barzilai-Borwein (BB) gradient descent method. While the use of gradient methods is a standard tool for the numerical approximation of optimal control laws, see \cite{borzibook}, the use of the BB method for large-scale agent-based models is relatively recent, \cite{ullmas}, and we report on its use as a reliable method for optimal consensus control problems (OCCP) in nonlinear MAS. While the control performance is satisfactory, our setting allows the controller to act differently on every agent at every instant. The question of a more parsimonious control design remains open. As an extension of the proposed methodology, we address the finite horizon OCCP with a non-smooth, sparsity-promoting control penalisation. This control synthesis is sparse, acting on a few agents over a finite time frame, however its numerical realisation is far more demanding due to the lack of smoothness in the cost functional. To circumvent this difficulty, we propose a numerical realisation of the control synthesis via metaheuristics related to particle swarm optimisation (PSO), and to nonlinear model predictive control (NMPC). Finally, based on the works \cite{CFRT,carrillj,FS13,BFRS15,ACFK17,ABCK15,AK18}, we discuss the resulting \textsl{mean-field optimal control problem}: that obtained as the number of agents $N$ tends to $\infty$ and the micro-state $(x_i(t),v_i(t))$ of the population is replaced by an agent density function $f(t,x,v)$. 

\textsl{Structure of the paper.} In Section 2 we revisit the Cucker-Smale model and results on consensus emergence and control. In Section 3 we address the OCCP  via first-order necessary optimality conditions and its numerical realisation. Section 4 introduces the sparse OCCP and its approximation. Concluding, in Section 5 we present a mean-field modelling approach for OCCP  when the number of agents is sufficiently large.  
 

\section{The Cucker-Smale model and consensus emergence}
We consider a set of $N$ agents with state $(x_i(t),v_i(t))\in\R^d\times\R^d$, where $d$ is the dimension of the physical space, interacting under second-order Cucker-Smale  dynamics
\begin{align}
\frac{dx_i}{dt}&=v_i\,,\quad \frac{dv_i}{dt}= \frac{1}{N}\sum_{j=1}^N a(\|x_i-x_j\|)(v_j-v_i)\label{eq:MAS1}\\
\bx(0)&=\bx_0\,,\quad \bv(0)=\bv_0\,,\label{eq:MAS2}
\end{align}
where $a(r)$ is a communication kernel of the type
\begin{align}
a(r)=\frac{K}{(1+r^2)^{\beta}}\,,\qquad\beta\geq0\,,K>0\,,
\end{align}
and we use the notation $\bx(t)=(x_1(t),\ldots,x_N(t))^t$, $\bv(t)=(v_1(t),\ldots,v_N(t))^t\in\R^{dN}$. Both $\|\cdot\|$ and $\|\cdot\|_2$ are indistinctly used for the $\ell_2$-norm, while $\|\cdot\|_1$ stands for the $\ell_1$-norm. 
We will focus on the study of consensus emergence, i.e. the convergence towards a configuration in which
\begin{align}\label{eq:vConsensus}
v_i=\bar{v}=\frac{1}{N}\sum_{j=1}^Nv_j\,\quad \forall i.
\end{align}
Note that although the interaction kernel $a(r)$ decays as the distance between agents increasing, the interconnection topology remains unchanged, i.e., all the agents interaction with all the agents in the swarm at all times. For a system of the type \eqref{eq:MAS1}-\eqref{eq:MAS2}, a consensus configuration will remain as such without any external intervention, and positions will evolve in planar formation. The emergence of consensus as a self-organisation phenomenon, either by a sufficiently cohesive initial configuration $(\bx_0,\bv_0)$ or a strong interaction $a(r)$, is a problem of interest in its own right.  To study consensus emergence, it is useful to define the form $B:\R^{dN}\times\R^{dN}\to\R$
\[
B(\mathbf{w},\mathbf{v})=\frac{1}{2N^2}\sum_{i,j=1}^N\|w_i-v_j\|^2\,.
\]
Note that for a population in consensus, $B(\bv,\bv)\equiv0$.  A solution $(\mathbf{x}(t),\mathbf{v}(t))$ to \eqref{eq:MAS1}-\eqref{eq:MAS2} tends to a consensus configuration if and only if
$$V(t) := B(\mathbf{v}(t),\mathbf{v}(t))\rightarrow 0 \quad \text{ as } \quad t\rightarrow+\infty.$$ 
Analogously, we define $X(t):=B(\bx(t),\bx(t))$. We briefly recall some well-known results on self-organised consensus emergence.

\begin{thm}{Unconditional consensus emergence (\cite{CS,carrillj}).} Given an interaction kernel $a(r)=\frac{K}{(1+r^2)^{\beta}}$ with $K>0$ and $0\leq\beta\leq\frac12$, the Cucker-Smale dynamics \eqref{eq:MAS1}-\eqref{eq:MAS2} convergence asymptotically to consensus, i.e. $V(0)\leq e^{-\lambda t} V(t)$, for $\lambda>0$.
\end{thm}

\begin{thm}{Conditional consensus emergence (\cite{HaLiu,HaHaKim}).} For $a(r)=\frac{K}{(1+r^2)^{\beta}}$ with $K>0$ and $\frac12\leq\beta$, if
$$\sqrt{V(0)}<\int\limits_{\sqrt{X(0)}}^{+\infty}a(2\sqrt{N}s)\,ds\,,$$
	
	then the Cucker-Smale dynamics \eqref{eq:MAS1}-\eqref{eq:MAS2} convergence asymptotically to consensus.
\end{thm}

In this work we are concerned with inducing consensus through the synthesis of an external forcing term
$\bu(t)=(u_1(t),\ldots,u_N(t))^t$
in the form 
\begin{align}
\frac{dx_i}{dt}&=v_i\,,\label{eq:MASC1}\\
\frac{dv_i}{dt}&= \frac{1}{N}\sum_{j=1}^N a(\|x_i-x_j\|)(v_j-v_i)+u_i(t)\,,\label{eq:MASC2}\\
x_i(0)&=x_0\,,\quad v_i(0)=v_0\,,\qquad i=1,\ldots,N\,,\label{eq:MASC3}
\end{align}

where the  control signals $u_i\in\U=\{u: \R_+\longrightarrow U\}$ and $U$ a compact subset of $\R^d$. 

%

\section{The optimal control problem and first-order necessary conditions}\label{sc:controlProblem}

In this section we entertain the problem of obtaining a centralised forcing term $\bu(t)$ which will either induce consensus on an initial configuration $(\bx_0,\bv_0)$ that would otherwise diverge, or which accelerates the rate of convergence for initial data that would naturally self-organise. Formally,  for $T>0$ and given a set of admissible control signals $\U^N: \R^+_0\to [L^\infty(0,T;\R^d)]^{N}$  for the entire population, we seek a solution to the minimisation problem
\begin{align}\label{eq:min}
\underset{\bu(\cdot)\in\U^N}{\min} \J(\bu(\cdot);\bx_0,\bv_0):=\int_0^T \ell(\bv(t),\bu(t))\,dt\,,
\end{align}
with the running cost defined as
\begin{align}
\ell(\bv,\bu):=\frac{1}{N}\sum_{j=1}^N\prt{\normsq{\bar{v}-v_j}+\gamma\normsq{u_j}}
,
\end{align}
with $\gamma>0$, subject to the dynamics \eqref{eq:MASC1}-\eqref{eq:MASC3}.

\subsection{First-order optimality conditions}
While existence of a minimiser $\bu^*$ of \eqref{eq:min} follows from the smoothness and convexity properties of the system dynamics and the cost, the Pontryagin Minimum Principle \cite{Pontryagin} yields first-order necessary conditions for the optimal control.  Let $(p_i(t),q_i(t))\in\R^d\times\R^d$ be adjoint variables associated to $(x_i,v_i)$, then the optimality system consists of a solution $(\bx^*,\bv^*,\bu^*,\mathbf{p}^*,\mathbf{q}^*)$ satisfying \eqref{eq:MASC1}-\eqref{eq:MASC3} along with the adjoint equations
\begin{align}
\label{eq:adjoint1}-\frac{dp_i}{dt}&=\frac{1}{N}\sum_{j=1}^N\frac{a'\prt{\norm{x_j-x_i}}}{\norm{x_j-x_i}}\inner{q_j-q_i}{v_j-v_i}\prt{x_j-x_i}\,,\\
\label{eq:adjoint2}-\frac{dq_i}{dt}&=p_i+\frac{1}{N}\sum_{j=1}^Na\prt{\norm{x_j-x_i}}\prt{q_j-q_i}-\frac{2}{N}\prt{\bar{v}-v_i},\\
\label{eq:adjoint3}p_i(T)&=0\,,\quad q_i(T)=0\,,\qquad i=1,\ldots,N\,,
\end{align}
and the optimality condition
\begin{align}\label{eq:optcond}
\bu(t)=\underset{\mathbf{w}\in\R^{dN}}{argmin}\sum_{j=1}^N\prt{\inner{q_j}{\frac{dv_j}{dt}}
+\frac{\gamma}{N}\normsq{w_j}}=-\frac{N}{2\gamma}\mathbf{q}^t\,.
\end{align}

\subsection{A gradient-based realisation of the optimality system}\label{sec:gradient}
The adjoint system \eqref{eq:adjoint1}-\eqref{eq:optcond} is used to implement a gradient descent method for the numerical realisation of the optimal control law. It can be readily verified that the gradient of the cost $\J$ in \eqref{eq:min} is given by
\begin{align}\label{eq:gradient}
\nabla\J\prt{\bu}=\mathbf{q}^t+\frac{2\gamma}{N} \bu\,,
\end{align}
obtained by differentiating \eqref{eq:optcond} with respect to $\bu$. With this expression for $\nabla\J$, the gradient iteration is presented in Algorithm \ref{alg:bb}.

\begin{algorithm}[H]
	\caption{Gradient descent method with BB update.}\label{alg:bb}
	\begin{algorithmic}
    \Require $tol>0$, $k_{\mathrm{max}}$, $\bu^0$, $\bu^{-1}$.
	\State $k=0$;
	\While{$\norm{\nabla\J\prt{\bu^{k}}}>tol$ and $k<k_{\mathrm{max}}$}
	\State 1) Obtain $(\bx^k,\bv^k)$ from \eqref{eq:MASC1}-\eqref{eq:MASC3} with $\bu^k$;
	\State 2) Obtain $(\mathbf{p}^k,\mathbf{q}^k)$ from \eqref{eq:adjoint1}-\eqref{eq:adjoint3} with $(\bx^k,\bv^k)$;
	\State 3) Evaluate the gradient $\nabla\J\prt{\bu^k}$ as in \eqref{eq:gradient};
	\State 4) Compute the step
	\begin{equation}\label{BBstep}
	\alpha_k := \frac{\inner{\bu^k-\bu^{k-1}}{\nabla\J\prt{\bu^k}-\nabla\J\prt{\bu^{k-1}}}}{\normsq{\nabla\J\prt{\bu^k}-\nabla\J\prt{\bu^{k-1}}}}\,;
	\end{equation}
	\State 5) Update $\bu^{k+1}=\bu^k-\alpha_k \nabla\J\prt{\bu^k}$;
	\State 6) $k:=k+1$;
	\EndWhile
	\end{algorithmic}
\end{algorithm}

In Algorithm 1, the gradient is first obtained by integrating the forward-backward optimality system, and then the step $\alpha_k$ in 4) is chosen as in the Barzilai-Borwein method, see \cite{Barzilai1988}.

\subsection{Numerical experiments}
\begin{center}
	\begin{figure}[h]
		\begin{center}
			\includegraphics[width=0.42\textwidth,height=.18\textwidth]{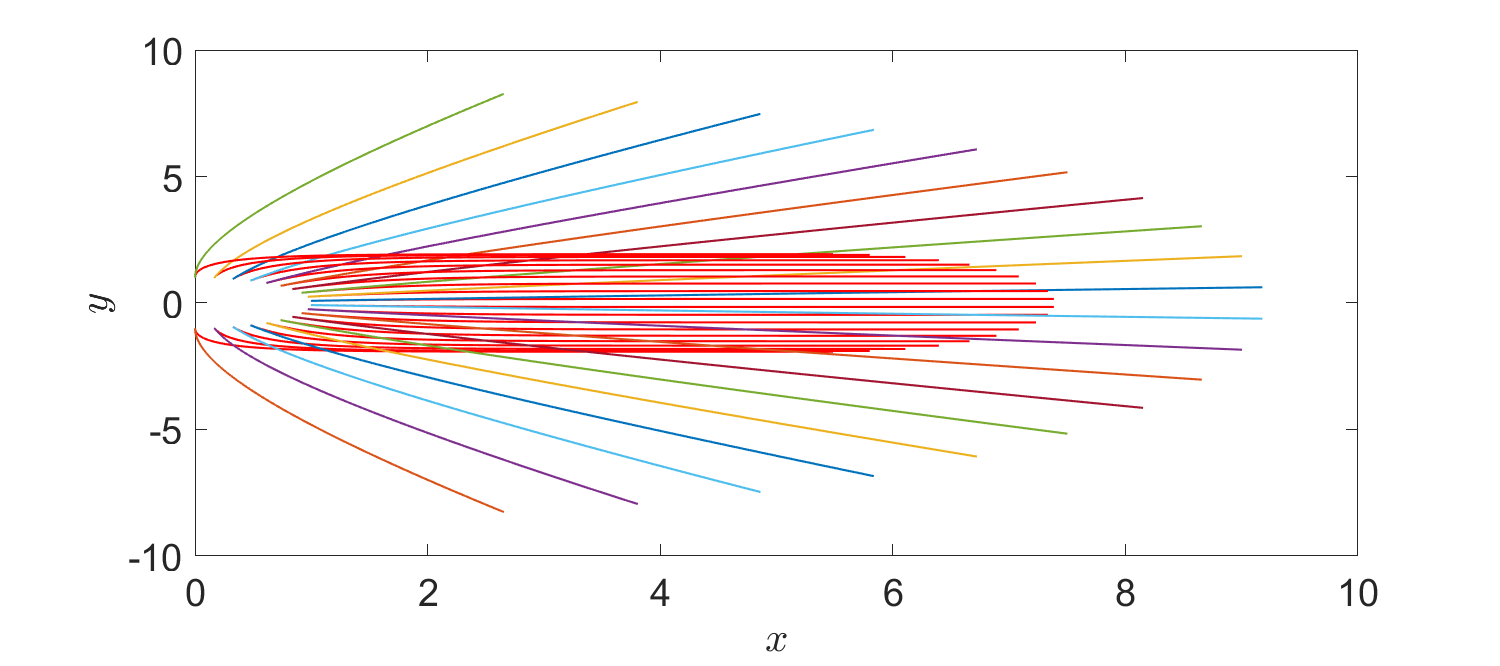}
			\includegraphics[width=0.42\textwidth,height=.18\textwidth]{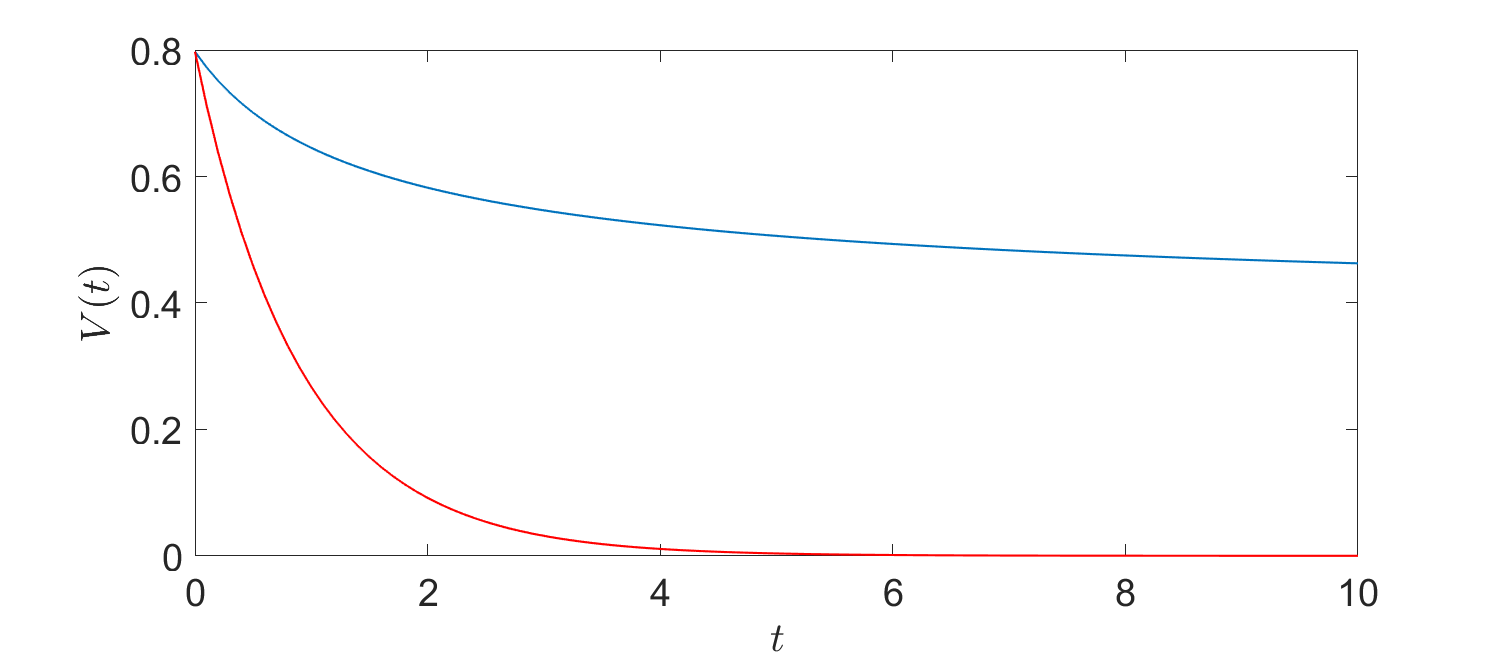}
			\includegraphics[width=0.42\textwidth,height=.18\textwidth]{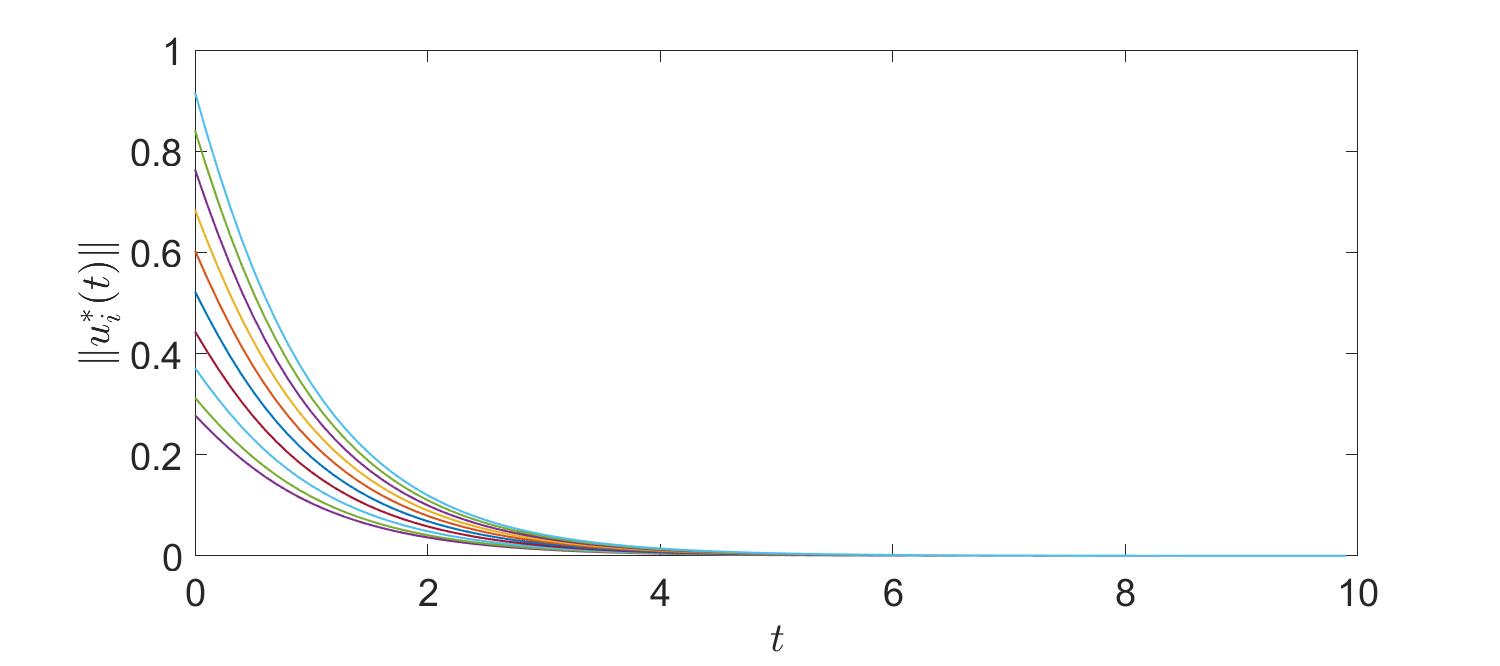}
			\caption{Free vs. controlled (shown in red) dynamics. Top: trajectories. Middle: functional $V(t)$. Bottom: norm of the optimal control $u_i^*(t)$ for each agent.}
			\label{fig:numericalSample}
		\end{center}
	\end{figure}
\end{center}

Figure \ref{fig:numericalSample} shows a comparison of the free two-dimensional dynamics of a sample initial condition $(\bx_0,\bv_0)$ and the system under an approximation to the optimal control $\bu^*$ found with this algorithm. A Runge-Kutta 4th order scheme was used to integrate the differential equations for the state and the adjoint with end time $T=10$, time step $dt=0.1$ (resulting on $N_T = 51$ points for the time discretisation), and a stopping tolerance for the gradient norm of $10^{-3}$. The condition $(\bx_0,\bv_0)$ is chosen such that consensus would not be reached naturally; long-time numerics of the free system show that $V(t)$ converges to an asymptotic value around $\bar{V}=0.4$. In the controlled setting, consensus is reached rapidly as can be seen from the trajectories themselves, as well as from the fast convergence of the functional $V(t)$ to zero. Furthermore, the norm of the control also decreases in time as the system is steered towards and into the self-organisation region.

\section{The sparse consensus control and its approximation via heuristics}

In this section, we address the problem of enforcing sparsity on the optimal consensus strategy. As shown in \cite{AFK17,BF,CFPT,KKR16}, one way to do it is by using as control cost the $\ell_1$-norm $\|\cdot\|_{1}$ in the minimisation problem \eqref{eq:min}, instead of the standard squared $\ell_2$-norm $\|\cdot\|_2^2$. However, the choice of the non-differentiable control cost $\|\cdot\|_{1}$ gives rise to a non-smooth cost functional $\mathcal{J}$, for which gradient-based numerical solvers like the one presented in Section \ref{sec:gradient} are not directly suitable. To circumvent the non-smoothness of $\mathcal{J}$, we shall resort to a metaheuristic procedure known as particle swarm optimisation (PSO).




\subsection{Particle Swarm optimisation}

First introduced in \cite{kennedy1995j,shi1998modified}, PSO is a numerical procedure that solves a minimisation problem by iteratively trying to improve a candidate solution. PSO solves the problem by generating a population of points in the discrete control state space of solutions $\mathcal{U}^N = \R^{d\times N\times N_T}$ called \textit{particles}. Each particle is treated as a point in this $D = d\times N\times N_T$-dimensional space with coordinates $(z_{i1},\ldots,z_{iD})$, and the cost functional is evaluated at each of these points. The best previous position (i.e., the one for which the cost functional is minimal) of any particle $(m_{i1},\ldots,m_{iD})$ is recorded, together with the index $h$ of the best particle among all the particles. We let then the particles evolve according to the system
\begin{align*}
z_{ij} &:= z_{ij} + w_{ij},\\
w_{ij} &:= w_{ij} + c_1\xi(m_{ij} - z_{ij}) + c_2\eta(m_{hj} - z_{ij}),
\end{align*}
where $c_1,c_2 > 0$ are two constant parameters and $\xi,\eta$ are two random variables with support in $[0,1]$. PSO is a metaheuristic since it makes few or no assumptions about the problem being optimised and can search very large spaces of candidate solutions: in particular, PSO does not use the gradient of the problem being optimised, which means PSO does not require that the optimisation problem be differentiable. However, it yields to a decrease of the cost function. Notice that whenever the dimension of the control space $D$ is very large (that is, either $d$, $N$ or $N_T$ are large), the problem suffers from the \textit{curse of dimensionality}, as the evaluation of $\mathcal{J}$ at all the particles and the subsequent search for the best position becomes prohibitively expensive. To mitigate this difficulty, we shall optimise within a nonlinear model  predictive control (NMPC) loop with short prediction horizon.

\subsection{Nonlinear Model Predictive Control}

For a prediction horizon of $H$ steps, for $k=1,\ldots,N_T-H$ and a discrete time version of the dynamics \eqref{eq:MAS1}-\eqref{eq:MAS2}, we minimise the following performance index
\begin{align}\label{eq:NMPC}
\sum^{H}_{h = 0}\frac{1}{N}\sum^N_{j = 1} \bigg( \|\overline{v}^{k+h} - v^{k+h}_j \|_2^2 + \gamma\|u^{k+h}_j\|^r_r \bigg),
\end{align}
for some $\ell_r$-norm with $r\geq 1$, generating  a  sequence  of  controls $(\mathbf{u}^k,\mathbf{u}^{k+1},\ldots,\mathbf{u}^{k+H})$ from which  only  the  first  term $\mathbf{u}^k$  is  taken  to  evolve  the  dynamics  from $k$ to $k+1$. The system' state is sampled again and the calculations are repeated starting from the new current state, yielding a new control and a new predicted state path. 
Although this approach is suboptimal with respect to the full time frame optimisation presented in Section 3, in practice it produces very satisfactory results.

For $H = 1$,  the NMPC approach recovers an instantaneous controller,  whereas for $H = N_T -1$ it  solves  the  full  time  frame  problem  \eqref{eq:min}.  Such  flexibility  is  complemented  with a robust behaviour,  as the optimisation is re-initialised every time step,  allowing to address perturbations along the optimal trajectory. For further references, see \cite{mayne2000constrained}.

\subsection{Numerical experiments} We now report the results of the numerical simulations of \eqref{eq:min} with $\|\cdot\|_1$ and $\|\cdot\|_2^2$ together with the setup PSO-NMPC described above. The aim is to check whether the optimal control obtained with the $\ell_1$-control cost is sparser than the one obtained with the $\ell_2$-norm. To do so, we shall compare their norms at each time, since a sparse control will be equal to 0 most of the time. Starting from an initial configuration $(\mathbf{x}_0,\mathbf{v}_0)$ that does not converge to consensus, we compare the effect of a different NMPC horizon $H \in \{3,10\}$ and of a different $\ell_r$-control cost for $r\in\{1,2\}$ in \eqref{eq:NMPC} on the optimal control strategy. 

\begin{center}
	\begin{figure}[h]
		\begin{center}
			\includegraphics[width=.45\textwidth,height=0.18\textwidth]{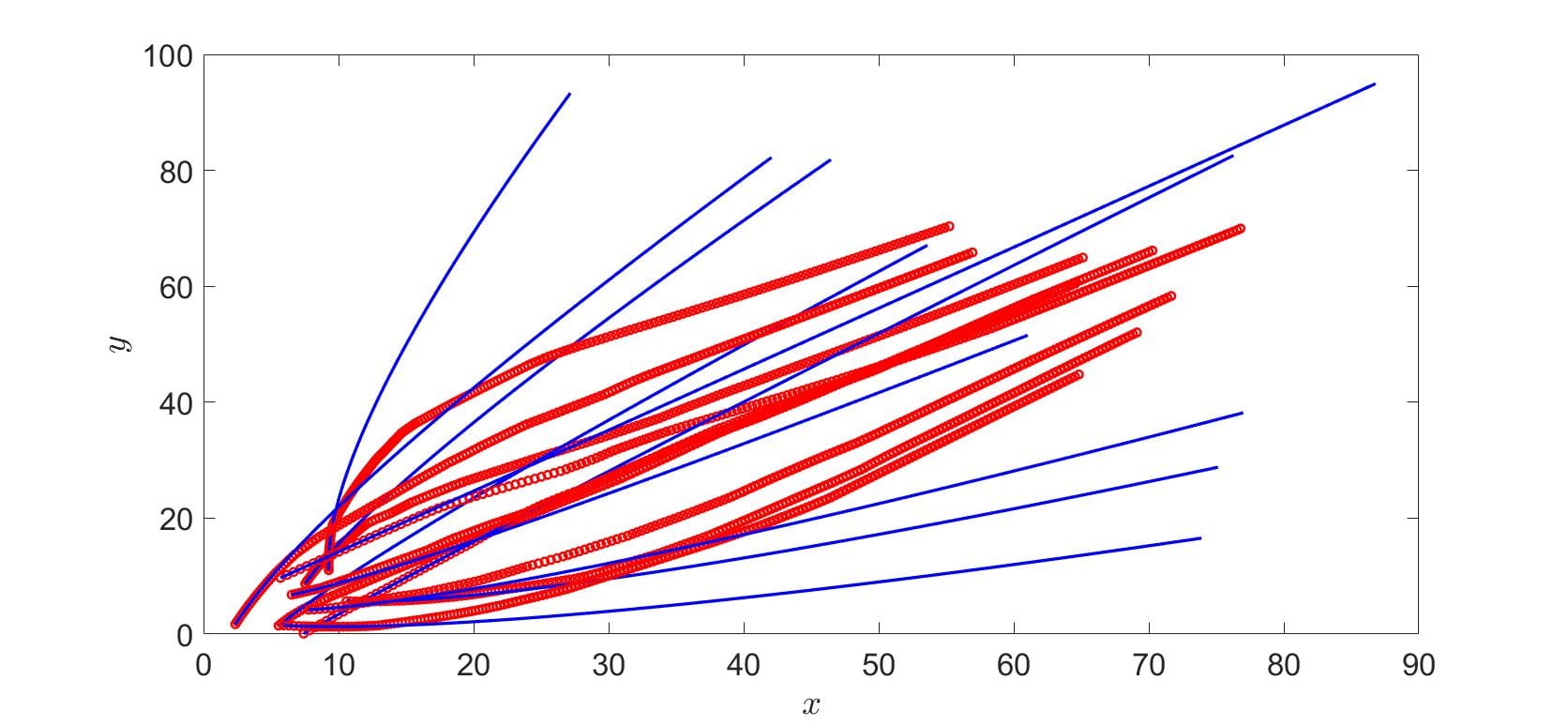}
			\includegraphics[width=.45\textwidth,height=0.18\textwidth]{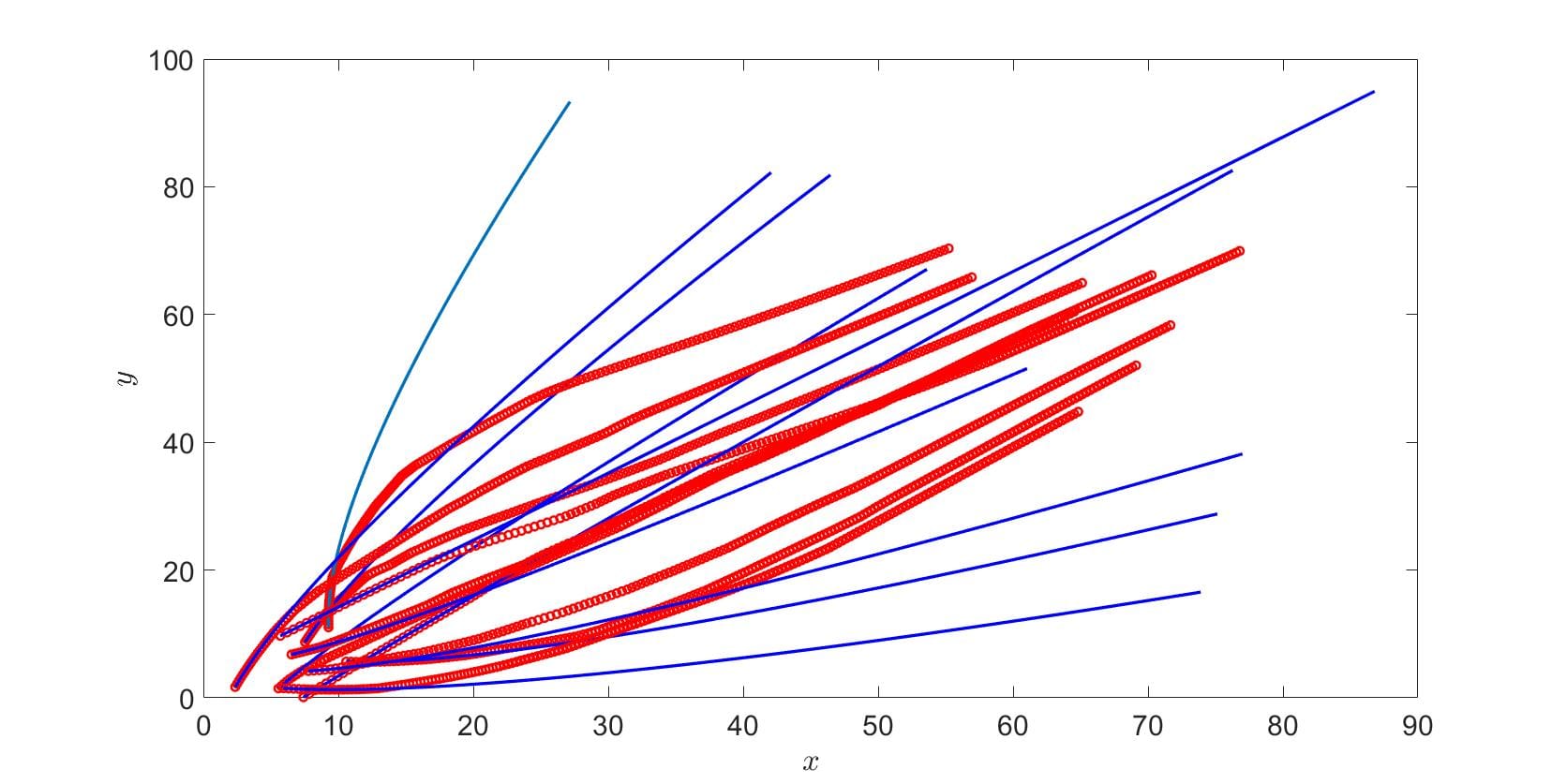}
			\caption{Controlled dynamics for $H = 3$. Top: agents' dynamics for $r=1$. Bottom: agents' dynamics for $r=2$.}
			\label{fig1}
		\end{center}
	\end{figure}
\end{center}
\begin{center}
	\begin{figure}[h]
		\begin{center}
			\includegraphics[width=0.45\textwidth,height=0.18\textwidth]{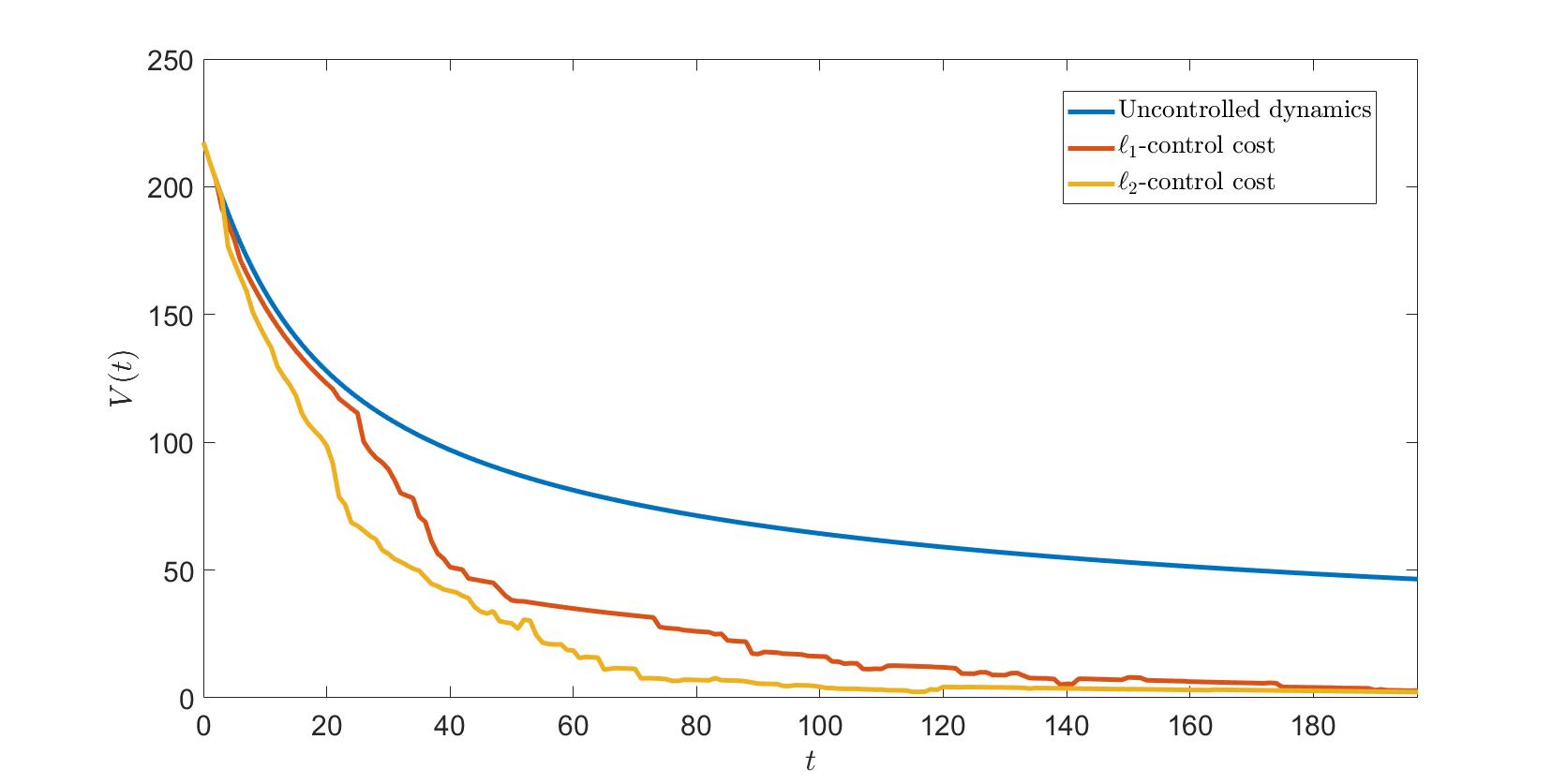}
			\caption{Functional $V(t)$ for $H = 3$. Blue: uncontrolled dynamics. Red: $\ell_1$-cost. Yellow: $\ell_2$-cost.}
			\label{fig2}
		\end{center}
	\end{figure}
\end{center}
\begin{center}
	\begin{figure}[h]
		\begin{center}
			\includegraphics[width=.44\textwidth,height=0.16\textwidth]{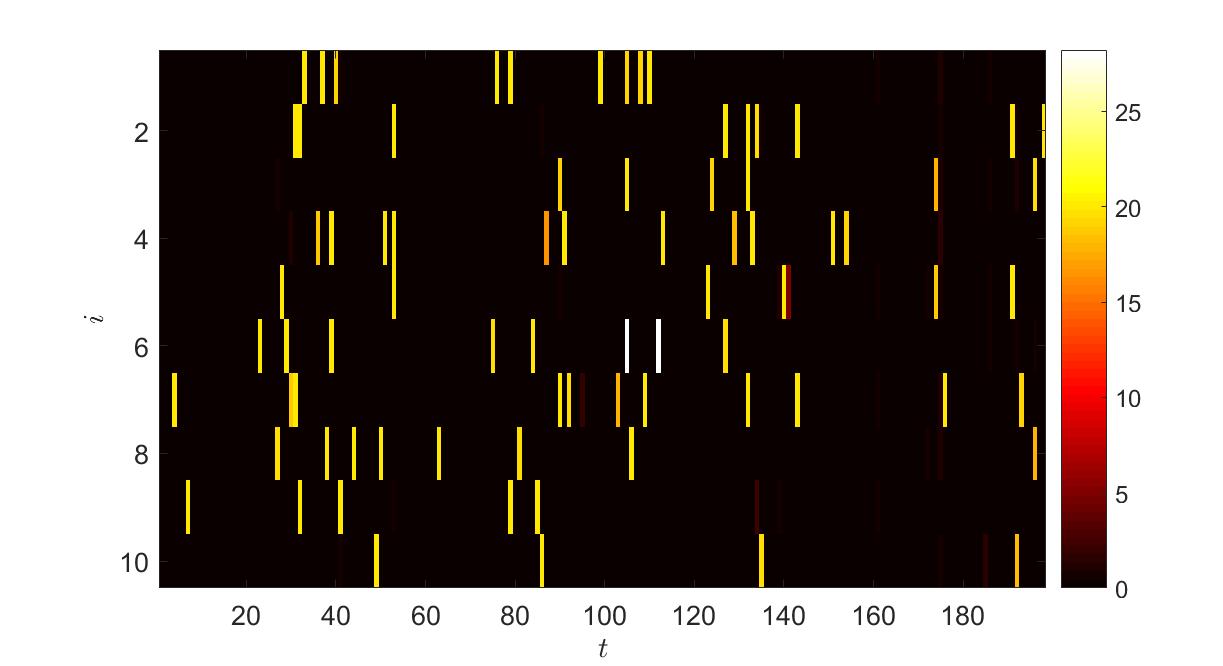}
			\includegraphics[width=.435\textwidth,height=0.16\textwidth]{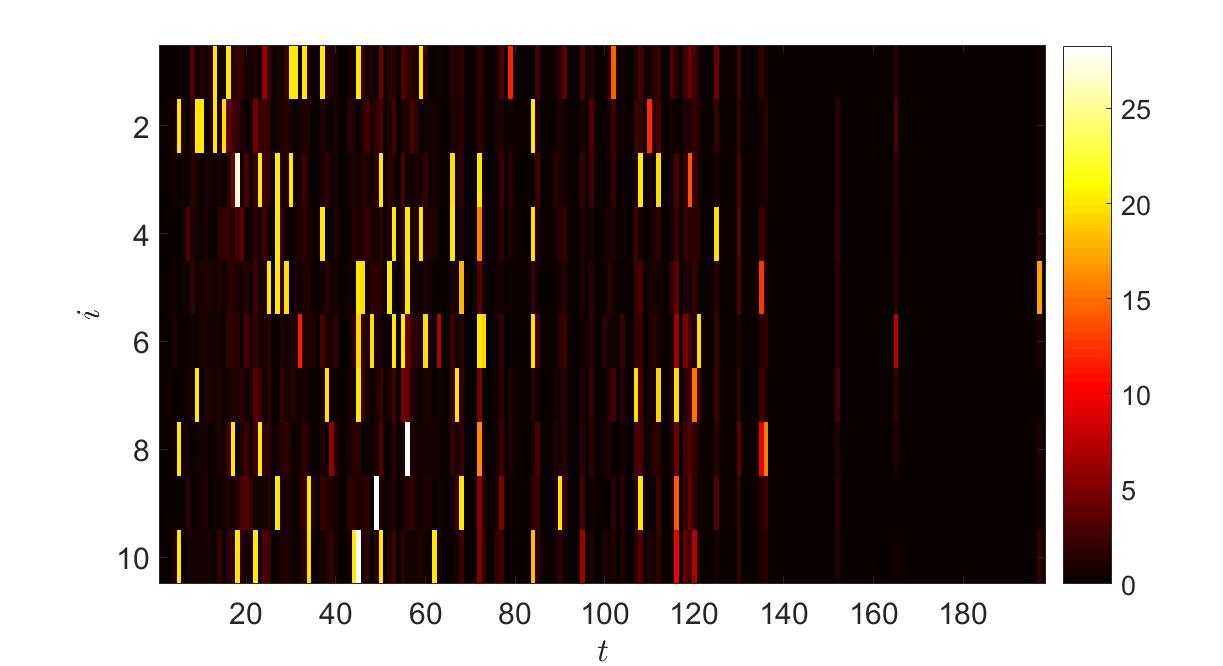}
			\caption{Heat map of control strategies for $H = 3$, control action in time for each agent. Top: the case of $r=1$. Bottom: the case of $r=2$.}
			\label{fig2bis}
		\end{center}
	\end{figure}
\end{center}

We test the PSO-NMPC procedure with $H=3$ periods ahead instead of the full time frame. Figure \ref{fig1} shows the controlled dynamics of the agents with the optimal control obtained for $r=1$ (top) and $r=2$ (bottom). Both controls decisively improve the alignment behaviour with respect to the uncontrolled dynamics. Figure \ref{fig2} shows the behaviour of the $V$ functional: for both controlled dynamics, the velocity spread $V(t)$ goes steadily to 0. To see how sparse the controls are, for each control strategy in Figure \ref{fig2bis} we show the corresponding \textit{heat map}, i.e., the matrix $\mathcal{H}$ such that $\mathcal{H}_{ik}$ contains the norm of the control acting on the $i$-th agent at step $k$. The stronger the control, the brighter the entry shall be: we can notice that the heat map for $r=1$ is sparser than the one for $r = 2$, being concentrated on few bright spots. This corroborates the findings of \cite{BF,CFPT}, where the sparsifying powers of an $\ell_1$-control cost were shown.

\section{The Continuous Control Problem}

\begin{center}
\begin{figure}[h]
	\begin{center}
		\includegraphics[width=0.42\textwidth,height=.18\textwidth]{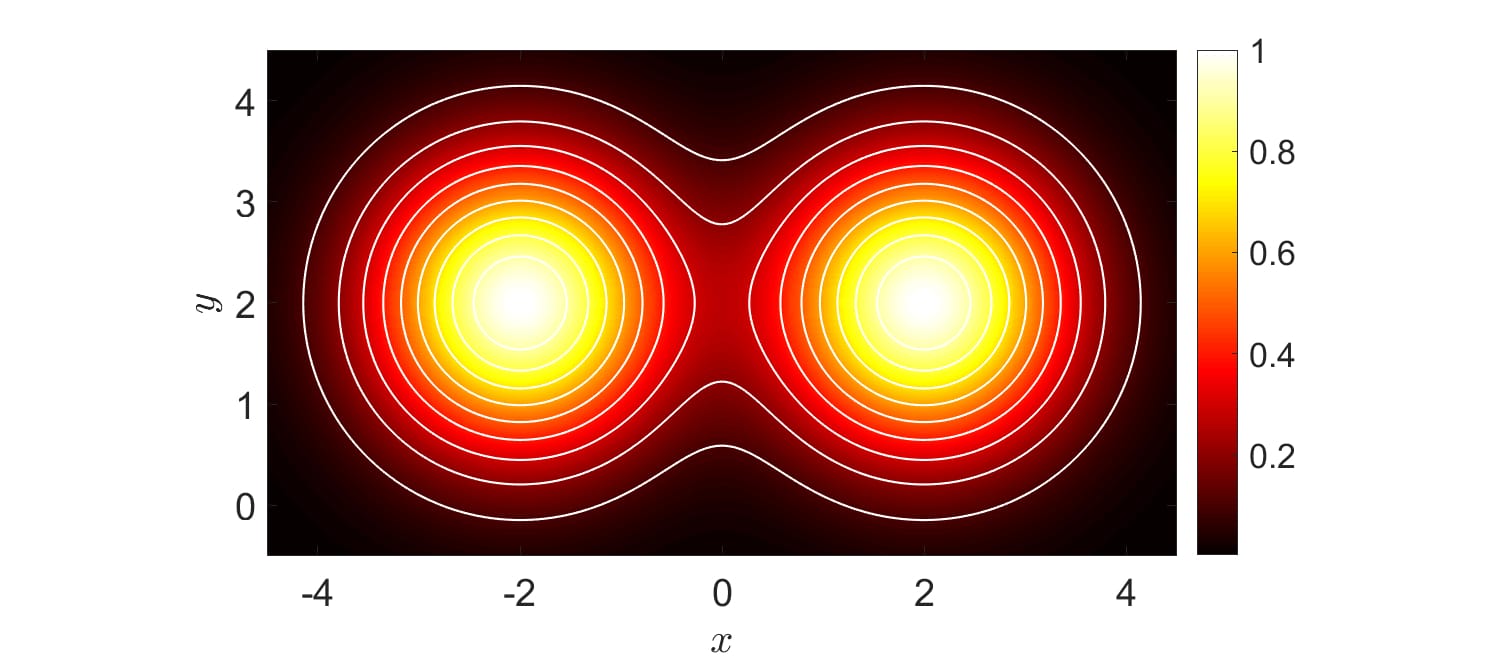}
		\caption{Marginal distribution of $f_0$ on $x$ and $v$.}
		\label{fig:initialCondition}
	\end{center}
\end{figure}
\end{center}

\begin{center}
\begin{figure}[h]
	\begin{center}
		\includegraphics[width=0.42\textwidth,height=.18\textwidth]{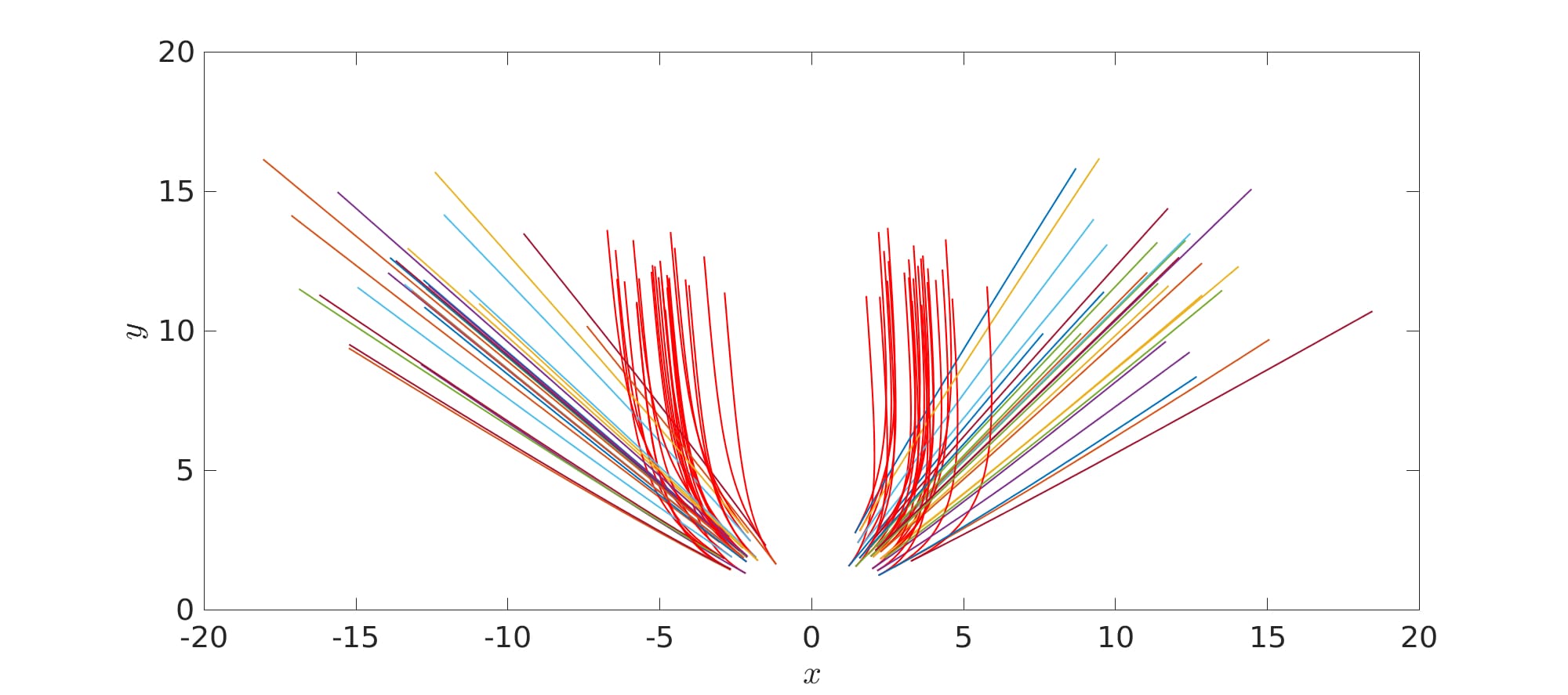}
		\includegraphics[width=0.42\textwidth,height=.18\textwidth]{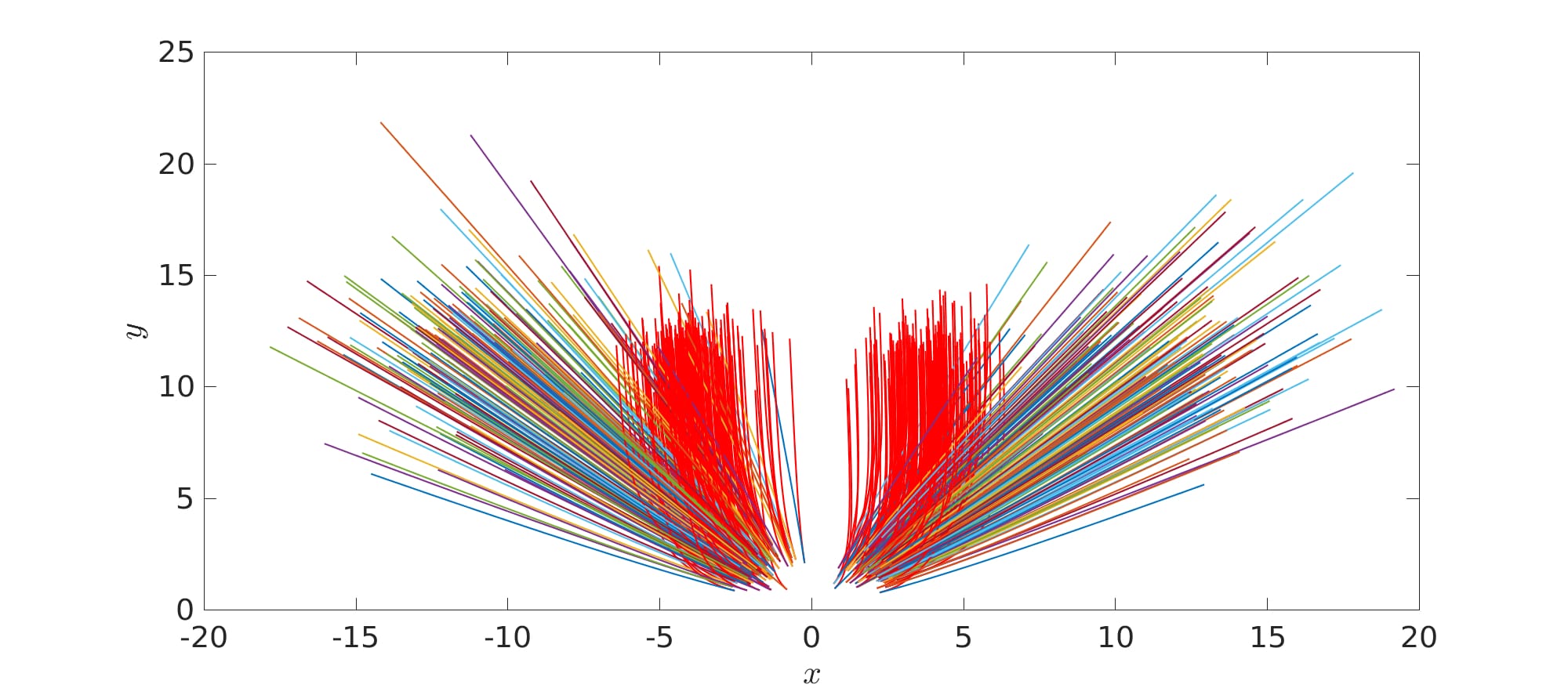}
		\includegraphics[width=0.42\textwidth,height=.18\textwidth]{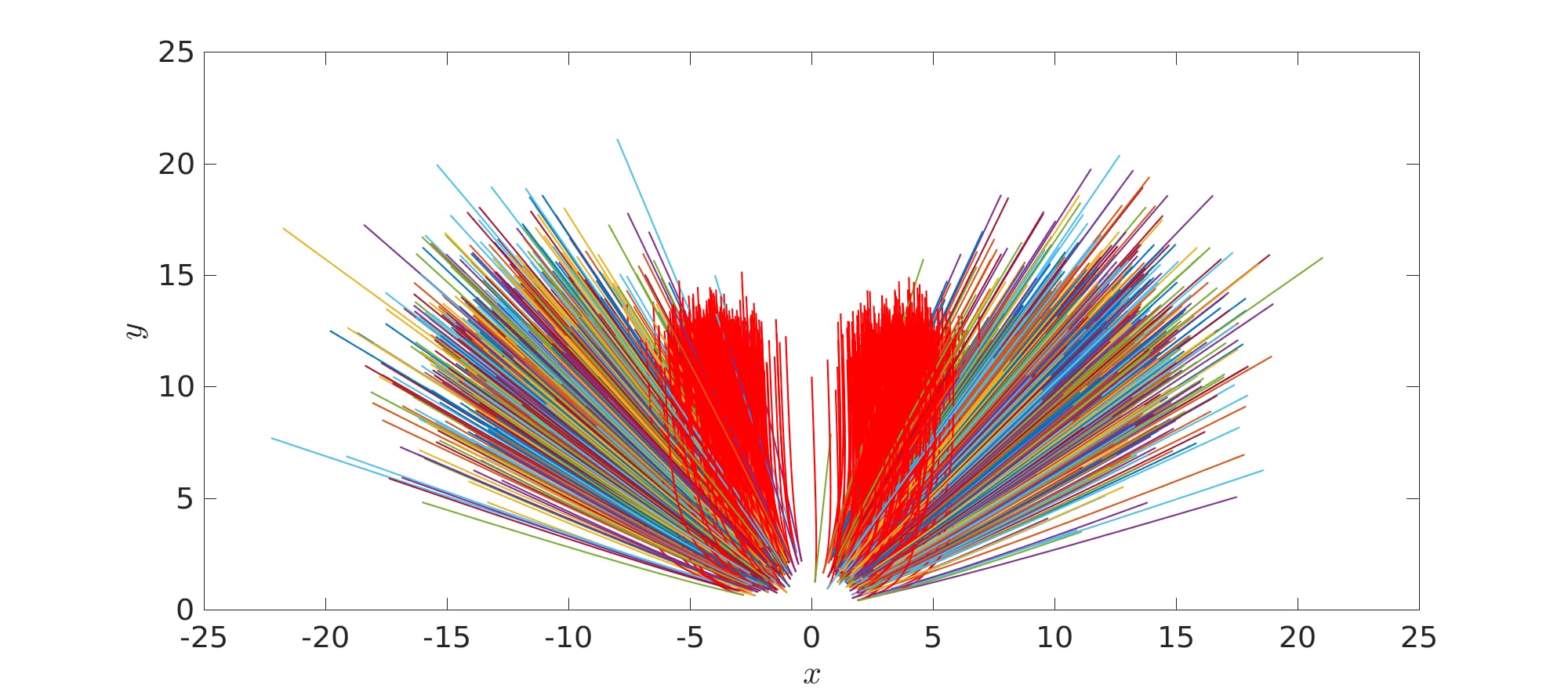}
		\caption{Free vs. optimally controlled (shown in red) dynamics. Top: $N=50$. Middle: $N=400$. Bottom: $N=2000$.}
		\label{fig:numericalMeanFieldSample}
	\end{center}
\end{figure}
\end{center}

\begin{center}
\begin{figure}[h]
	\begin{center}
		\includegraphics[width=0.42\textwidth,height=.18\textwidth]{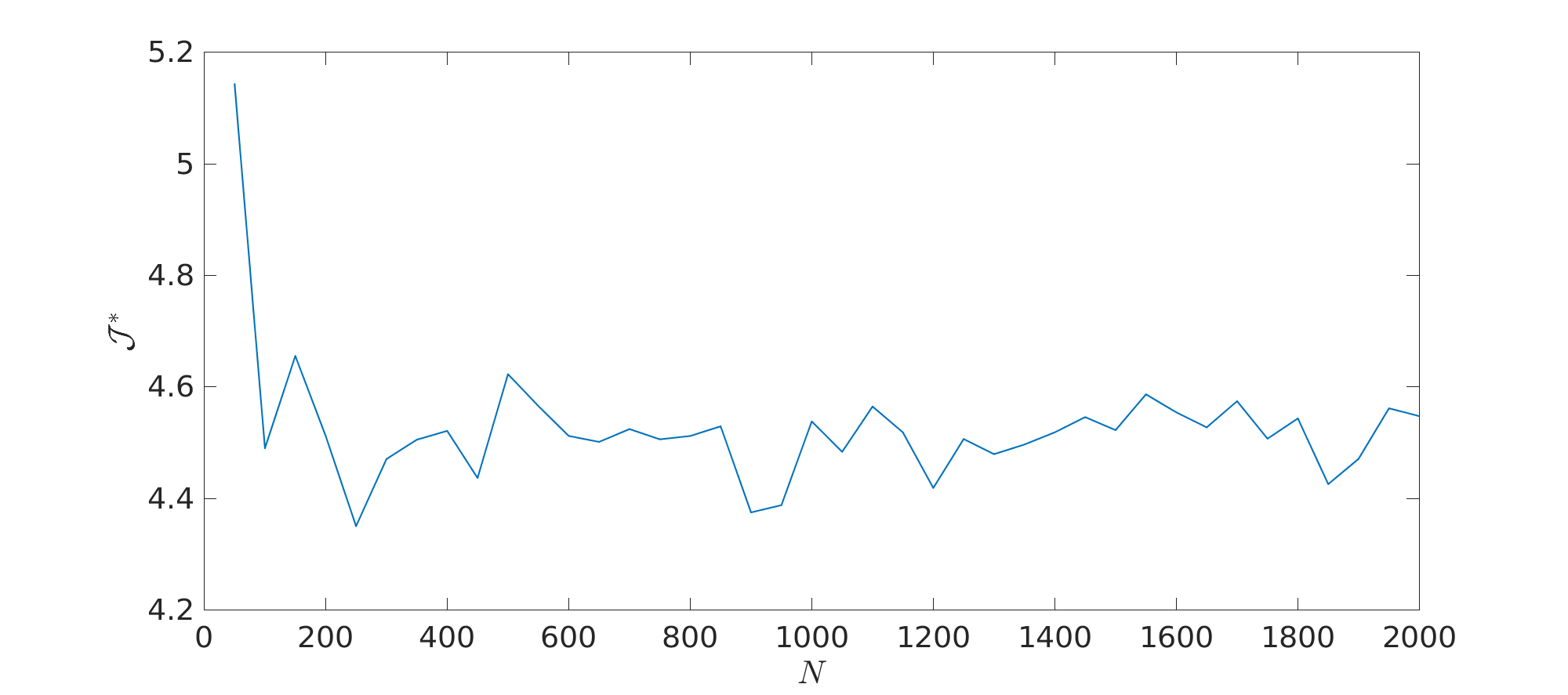}
		\caption{Evolution of $\J^*$ with $N\in\left\{50,100,\cdots,2000\right\}$.}
		\label{fig:meanFieldU}
	\end{center}
\end{figure}
\end{center}

\begin{center}
\begin{figure}[h]
	\begin{center}
		\includegraphics[width=0.42\textwidth,height=.18\textwidth]{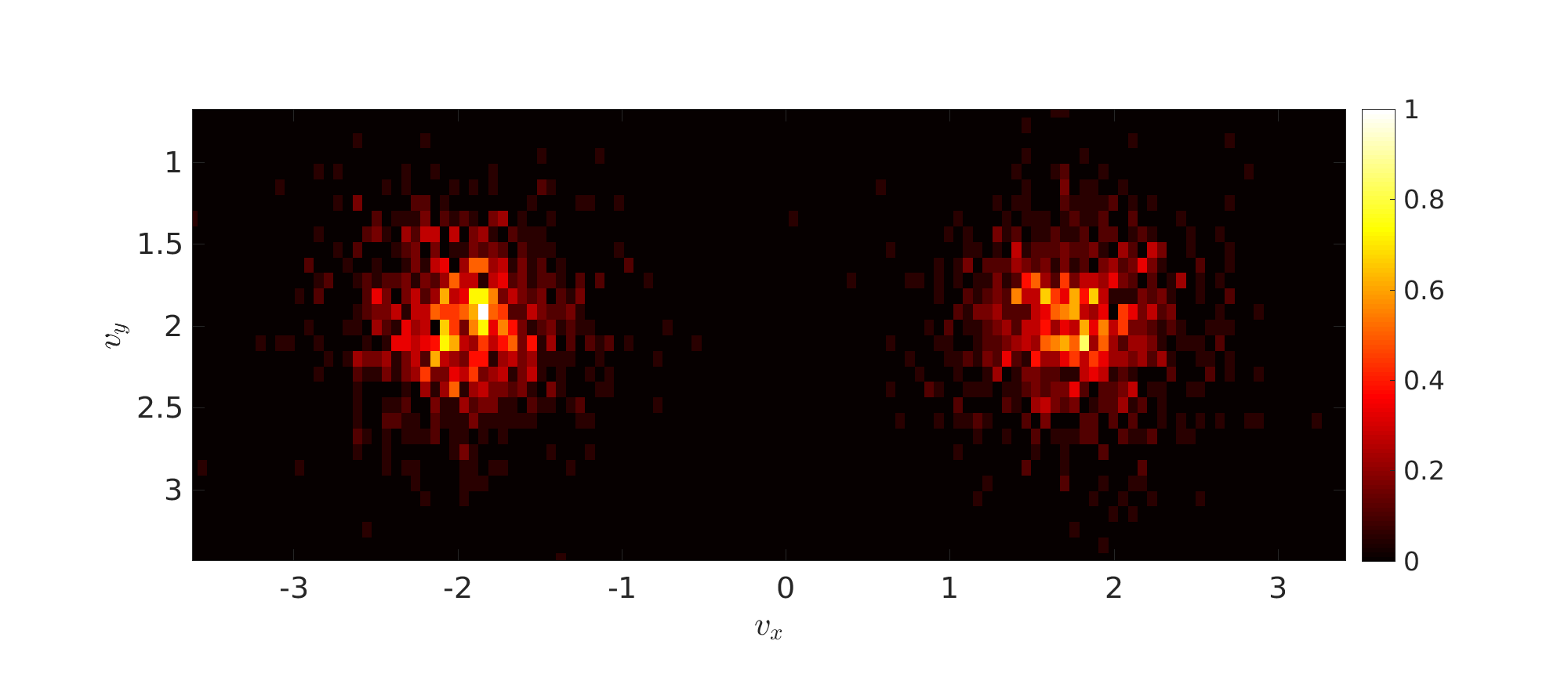}
		\includegraphics[width=0.42\textwidth,height=.18\textwidth]{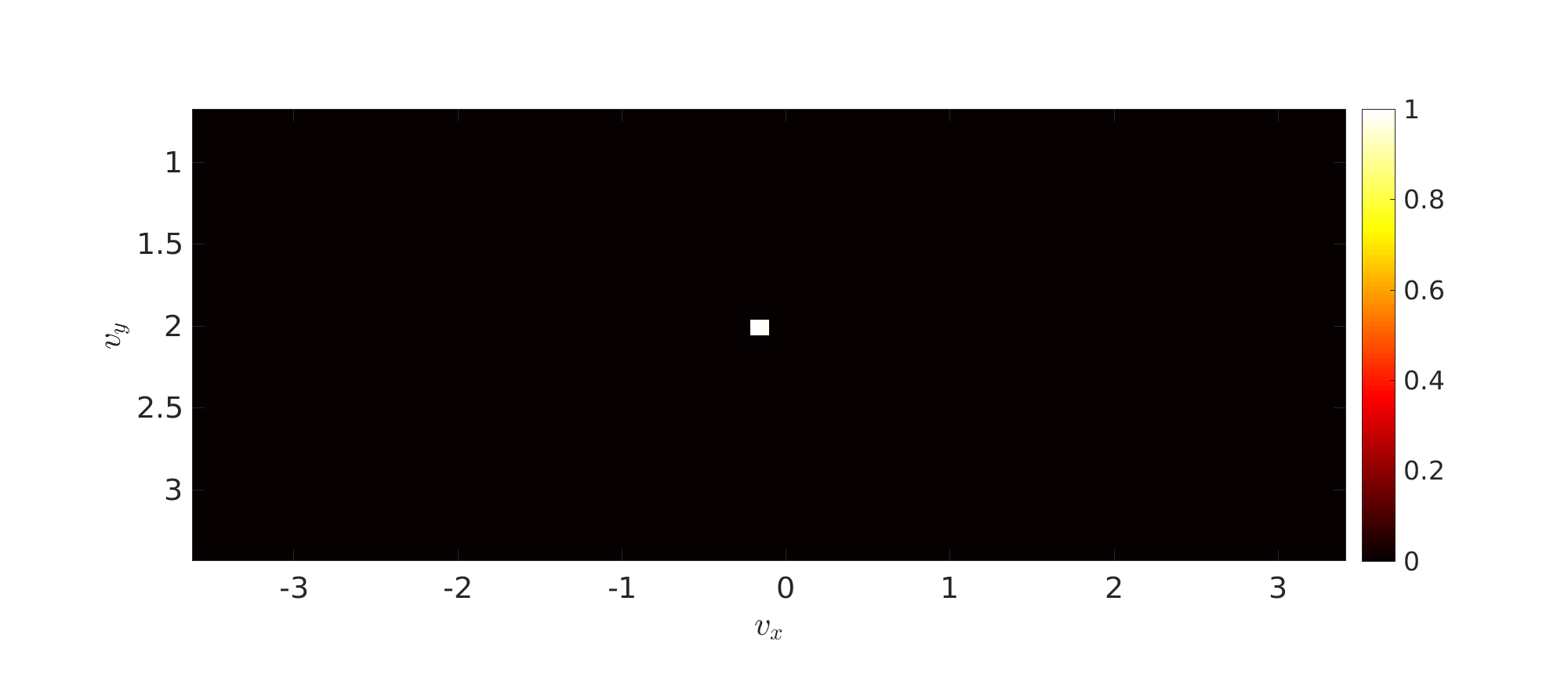}
		\caption{Marginal distributions of $f(T)$ on $v$, $N=2000$. Top: free setting. Bottom: controlled setting.
			Consensus emergence produces a concentration of the distribution around a single velocity close to $(0,2)$.
		}
		\label{fig:meanFieldF}
	\end{center}
\end{figure}
\end{center}

\begin{center}
\begin{figure}[h]
	\begin{center}
		\includegraphics[width=0.42\textwidth,height=.18\textwidth]{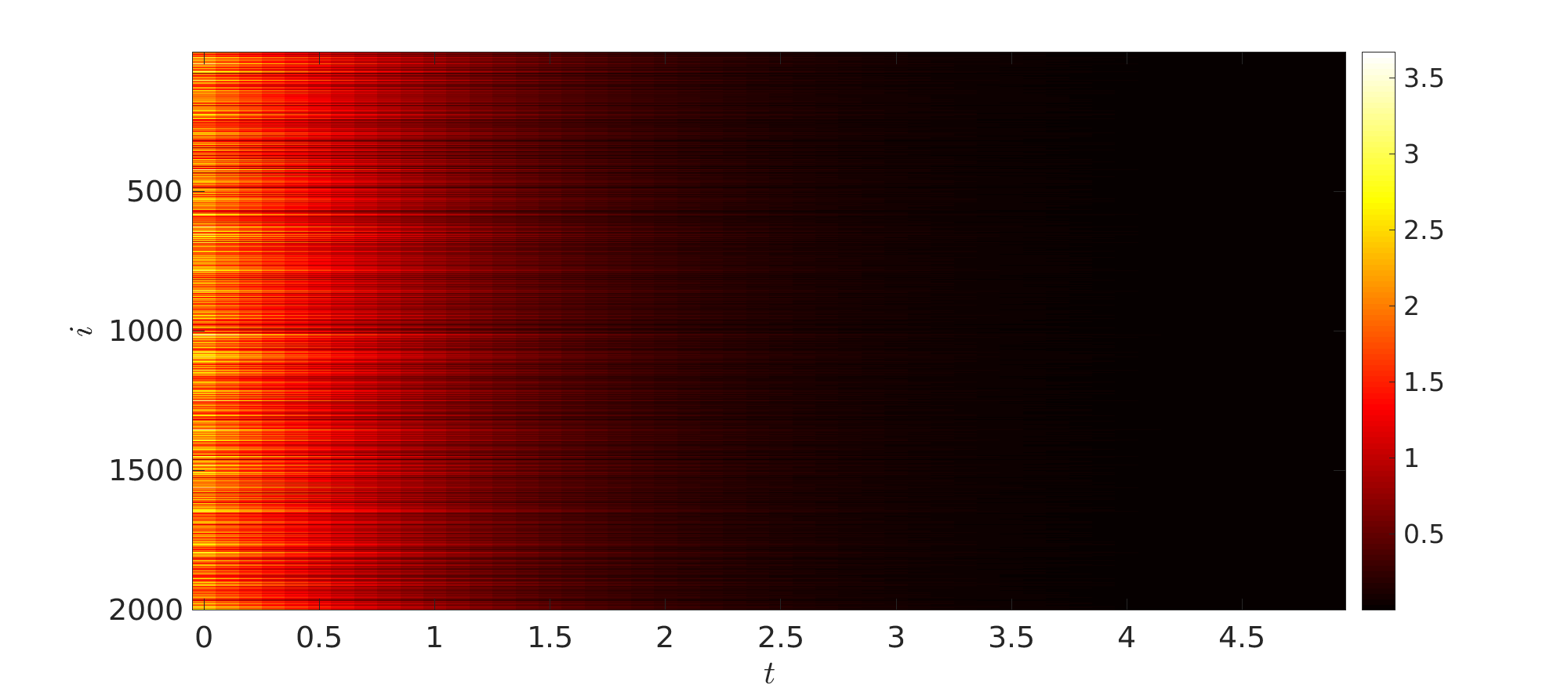}
		\caption{Heat map of $\norm{u_i^*(t)}$ with $N=2000$.}
		\label{fig:meanFieldUHeat}
	\end{center}
\end{figure}
\end{center}

We consider the continuous control problem that results in the limit of the discrete problem from Section \ref{sc:controlProblem} as $N\rightarrow\infty$. Formally, the forced Cucker-Smale dynamics \eqref{eq:MASC1}-\eqref{eq:MASC3} can be written as a Vlasov-type transport equation
\begin{equation}\label{eq:meanfield}
\frac{\partial f}{\partial t}+\prt{v\cdot\nabla_x}f+\nabla_v\cdot\brk{\prt{A(x,v)*f+u(t,x,v)}f}=0,
\end{equation}
where $A(x,v)=a(|x|)v$ and $f, u:[0,T]\times\R^2\times\R^2\rightarrow \R$ are the probability density for the state and a forcing term, respectively. An equivalent minimisation problem can be posed:
\begin{align}\label{eq:meanFieldMin}
\underset{u\in\U}{\min} \J(u;f_0):=\int_0^T \ell(f,u)\,dt\,,
\end{align}
for fixed $T>0$ and a running cost $\ell(f,u)$ given by the expression:
\begin{align}
\int_{\R^{2d}}
\normsq{v-\int_{\R^{2d}}w\,df(y,w)}+
\gamma\normsq{u(t,x,v)}
\,df(x,v).
\end{align}
The solutions of the discrete control problem \eqref{eq:min} converge to that of the continuous problem \eqref{eq:meanfield}, as discussed in \cite{FS13}. This can be verified numerically by fixing an initial distribution $f_0$ and studying the sequence solutions of the discrete problem with initial conditions sampled from said distribution; a subsequence is known to converge as $N\rightarrow\infty$. Besides the solution, the optimal value of the objective functional $\J^*_N$ is also expected to converge, which can be verified. Initial conditions without natural consensus were constructed by sampling $x$ from a superposition of two Gaussian distributions and letting $x=v$. The marginal distributions of $f_0$ on $x$ and in $v$ are shown in Figure \ref{fig:initialCondition}.
A sequence of such discrete problems were solved for several values of $N$. Figure \ref{fig:numericalMeanFieldSample} shows the comparison between the free and controlled trajectories for various values of $N$. The Runge-Kutta scheme was used to solve the differential equations for the state and the adjoint with end time $T=5$, time step $dt=0.1$ (resulting on $51$ points for the time discretisation), and a stopping tolerance of $10^{-2}$. Figure \ref{fig:meanFieldU} shows the evolution of the optimal cost $\J^*_N$, which appears to be of order $\J^*_N\sim \mathcal{O}(1)$ as expected for the convergence as $N\rightarrow\infty$. Figure \ref{fig:meanFieldF} shows the marginal distribution of $f(T)$ on $v$ for the free and forced settings with the same scale; the controlled case yields a singular distribution indicating consensus. Figure \ref{fig:meanFieldUHeat} shows a heat map of the optimal control $u_i^*(t)$. We observed that the average norm of the control is of $\sim \mathcal{O}(1)$ as $N\rightarrow\infty$; furthermore the time at which the control is nearly zero is roughly constant for large $N$. Table \ref{table:computationTime} shows the evolution of the number of optimisation iterations (i.e. loops on Algorithm \ref{alg:bb}) as well as the computation CPU time in hours; notice that the number of iterations remains roughly constant, while the computation time scales quadratically in $N$.


\input{meanFieldTable.tex}

\section{Acknowledgements}
JAC was partially supported by the EPSRC grant EP/P031587/1.

\bibliography{biblioACFK}     
       
\end{document}

%% file: meanFieldTable.tex
\begin{table}
\centering
\begin{tabular}{cccccccccccc}
Agents ($\mathbf{N}$)
& 50   & 100  & 150  & 200  & 250  & 300 \\
Iterations ($\mathbf{I}$)
& 26 & 25 & 24 & 27 & 28 & 28 \\
Time ($\mathbf{T}$)
& 0.1 & 0.1 & 0.3 & 0.7 & 1.1 & 1.4 \\


%
%

&&&&&&\\$N$
& 1550 & 1600 & 1650 & 1700 & 1750 & 1800 \\$I$
& 28 & 28 & 28 & 28 & 28 & 28 \\$T$
& 49.1 & 51.7 & 52.8 & 58.8 & 59.8 & 65.0 \\

&&&&&&\\$N$
& 1850 & 1900 & 1950 & 2000 \\$I$
& 27 & 29 & 28 & 28 \\$T$
& 65.4 & 75.3 & 77.4 & 82.8 \\

&&&&&&\\
\end{tabular}
\caption{Number of iterations ($I$) and CPU time ($T$) in hours for $N$ agents for the discrete optimisation problem.}
\label{table:computationTime}
\end{table}